\documentclass[12pt]{amsart}
\usepackage{epsfig,color}

\makeatother

\theoremstyle{definition}

\def\fnum{equation}

\numberwithin{equation}{section}
\newcommand{\Vol}{{\text{Vol}}}
\newcommand{\V}{{\text{V}}}

\newcommand{\Ric}{{\text{Ric}}}

\newcommand{\Hess}{{\text {Hess}}}

\def\RR{{\bold R}}

\newcommand{\e}{{\text {e}}}

\newcommand{\eqr}[1]{(\ref{#1})}

\title[Monotonicity - analytic and geometric implications]{Monotonicity - analytic and geometric implications}

 \author{Tobias Holck Colding}%
\address{MIT, Dept. of Math.\\
77 Massachusetts Avenue, Cambridge, MA 02139-4307.}

\author{William P. Minicozzi II}%
\address{Johns Hopkins University\\
Dept. of Math.\\
3400 N. Charles St.\\
Baltimore, MD 21218.}

\thanks{The authors
were partially supported by NSF Grants DMS  11040934, DMS
0906233,  and NSF FRG grants DMS 
 0854774 and DMS 0853501.  This material is based upon work supported by the National Science 
 Foundation under Grant No. 0932078, while the first author was in residence at 
 the Mathematical Science Research Institute (MSRI) in Berkeley, 
 California, during the Fall of 2011.}

\email{colding@math.mit.edu and minicozz@math.jhu.edu}

\begin{document}

\maketitle

\begin{abstract} 
In this expository article, we discuss various monotonicity formulas for parabolic and elliptic operators and explain how the analysis of function spaces and the geometry of the underlining spaces are intertwined.  

\hskip4mm 
After briefly discussing some of the well-known analytical applications of monotonicity for parabolic operators, we turn to their elliptic counterparts, their geometric meaning, and some geometric consequences. 
 \end{abstract}

\setcounter{section}{-1}

\section{Introduction}

The aim of this survey is to explain some new and old monotonicity formulas and describe some of their applications.  These formulas can roughly be divided into three groups: elliptic, parabolic and geometric.   Each has their own advantages and  applications and we will describe some of those.   For instance, we will touch upon how the parabolic monotonicity formulas imply functional inequalities and the elliptic ones imply uniqueness of blow-ups,
 whereas the geometric ones imply compactness for the spaces in question and cone structure for those same spaces.  But our real focus is how these different monotonicity formulas are linked and all have both analytic and geometric implications. 

\section{Functional inequalities and monotonicity}

In this section, we discuss some sharp monotonicity formulas and sharp gradient estimates for the heat and Laplace equations on manifolds and mention briefly some of the applications of the parabolic formulas.  The next section will focus on the interplay between them and the local geometry and, in particular, the question of uniqueness of blow-ups or blow-downs of the space.

\subsection{Parabolic operators}

If $f>0$ is a positive function on a Riemannian manifold $M$ with $\int_Mf=1$, then the Shannon entropy $S_0$
 and Fisher information $F_0$ are defined as follows
\begin{align}
S_0&=-\int_M \log f\, f\, , \notag\\
F_0&= \int_M \frac{|\nabla f|^2}{f}\, .\notag
\end{align}
Similarly,  if $u>0$ is a solution to the heat equation $(\partial_t-\Delta)\, u=0$ with $\int_Mu=1$, then we define $S(t)$ by
\begin{align}
S(t) &=S_0(t)-\frac{n}{2} \log (4\pi t)-\frac{n}{2} \notag \\
&=-\int_M \log u\, u-\frac{n}{2} \log (4\pi t)-\frac{n}{2}\, .
\end{align}
Here, $S$ is normalized so that   $S$ remains identically zero for all time when $u$
 is the heat kernel 
 \begin{align}
 	H(x,y,t)=(4\pi t)^{-\frac{n}{2}} \exp \left(-\frac{|x-y|^2}{4t}\right) \notag 
\end{align}
on Euclidean space $\RR^n$.   Our parabolic 
Fisher information $F$ will also differ from $F_0$ by a normalization that comes from the Euclidean heat kernel.  
Taking the derivative of the entropy $S$, using the heat equation, and integrating by parts gives $1/t$ times $F$, where $F$ is  
\begin{align}
F(t)=t\,F_0(t)-\frac{n}{2}=t \int_M \frac{|\nabla u|^2}{u}-\frac{n}{2}\, .
\end{align}
(We   multiplied   $F_0$ by $t$ so that $S$ and $F$ scale the same way.)

Following \cite{P2}, we set $W=S+F$ and get that on a Ricci-flat manifold or, more generally, on a manifold with $\Ric\geq 0$, that $W'\leq 0$\footnote{A space-time version of $W$ plays an important role in Perelman's work on the Ricci flow.}.  In fact, a computation shows that 
\begin{align}  \label{e:deriW}
W' &=\frac{(t\,F)'}{t}  \\
&= -2t\int_M\left( \left|\Hess_{\log u}+\frac{1}{2t}g\right|^2+\Ric (\nabla \log u,\nabla \log u)\right)\,u\, ; \notag
\end{align}
where $g$ is the Riemannian metric and $\Hess_{\log u}$ is the Hessian of $\log u$.  When $\Ric\geq 0$ this implies (as $t\,F$ is obviously $0$ for $t=0$) that $F\leq 0$ and hence $S\downarrow$.

We will next briefly see how one can use these quantities and formulas to prove some functional inequalities using monotonicity.  As examples, we single out the Cauchy-Schwarz inequality and the logarithmic Sobolev inequality.  Many other interesting inequalities can be proven by monotonicity; see, for instance, \cite{BL}, \cite{BCeLM}, \cite{BB}, \cite{BCarChrT}, \cite{Car}, \cite{CaLL}, \cite{CeL}, \cite{V}.  The  Cauchy-Schwarz inequality  is perhaps not only the simplest example of an inequality that can be proven from monotonicity using the heat equation, but the argument itself is also perhaps the easiest.  Another well-known instance is the log-Sobolev inequality.  

The log-Sobolev inequality of \cite{G} is that  for any function $f>0$ on $\RR^n$ with $(2\pi)^{-\frac{n}{2}} \int f^2\,\e^{-\frac{|x|^2}{2}}=1$, 
\begin{align}  \label{e:logsob0}
\int f^2\log f\,\e^{-\frac{|x|^2}{2}}\leq \int |\nabla f|^2\,\e^{-\frac{|x|^2}{2}}\, .
\end{align}
(The point being that the gain in the exponent in the usual Sobolev inequality is dimension dependent;   in the log-Sobolev inequality there is a gain of a log factor independent of the dimension.  There is a second point and that is that as $n\to \infty$ the measures $dx_1\wedge\cdots \wedge dx_n$ do not converge to a Radon measure whereas $(2\pi)^{-\frac{n}{2}} \e^{-\frac{|x|^2}{2}}\,dx_1\wedge\cdots \wedge dx_n$ do.)  By a change of variables (see, for instance, (2.1) on page 693 of 
\cite{BL}; cf. also \cite{Ca}),
 the log-Sobolev inequality is equivalent to that  for any function $f>0$ with $\int f=1$
\begin{align}  \label{e:logsob1}
\frac{n}{2} \log \left( \frac{1}{2n\pi\e}\int \frac{|\nabla f|^2}{ f}\right) &-\int f\,\log f \notag \\
	&=\frac{n}{2} \log \left( \frac{F_0}{2n\pi\e }\right)+S_0\geq 0\, ;
\end{align}
with equality for any Gaussian.  To show this inequality on $\RR^n$, let more generally $M$ be a manifold with $\Ric\geq 0$ and suppose that $u$ solves the heat equation with initial condition $u(\cdot,0)=f$.  The above inequality will follow from showing that for all  $t$
\begin{align}  \label{e:logsob1a}
\frac{n}{2} \log \left( \frac{F+\frac{n}{2}}{2n\pi\e }\right)+S+\frac{n}{2} \log (4\pi )+\frac{n}{2}\geq 0\, .
\end{align}
To show this, take the derivative of the left hand side of \eqr{e:logsob1a} to get
\begin{align}
\frac{n}{2t\,(F+\frac{n}{2})}\left((tF)'+\frac{2}{n}\, F^2\right)\leq 0\, .
\end{align}
Here, the inequality follows from the formula for $W'$ together with the Cauchy-Schwarz inequality; see, for instance, corollary 5.15 of \cite{C2}.
Briefly, to get the log-Sobolev inequality on $\RR^n$, use the monotonicity together with the asymptotic formula for solutions of the heat equation 
 (solutions on $\RR^n$ converge, after rescaling, at infinity to constant multiples of the usual Gaussian).  It is easy to see (more or less from the argument given above) that 
\eqr{e:logsob1} does not hold on a general manifold $M$ with $\Ric\geq 0$.  Rather, by \cite{BE} (see also \cite{B}),   if $f>0$ and $\int_Mf^2\,\e^{-\phi}=1$, then
\begin{align}  \label{e:logsob2}
\int_Mf^2\, \log f\, \e^{-\phi}\leq \frac{1}{c}\int_M |\nabla f|^2\, \e^{-\phi}\, , 
\end{align}
as long as $\Ric+\Hess_{\phi}\geq c\, g$ and $\int_M\e^{-\phi}=1$.

Another important and, as it turned out, closely related result is the Li-Yau gradient estimate, \cite{LY}.  This asserts that 
if $u>0$ is a solution to the heat equation $(\partial_t-\Delta)\, u=0$
on a manifold with nonnegative Ricci curvature,  then
\begin{align}
t\,\left(\frac{|\nabla u|^2}{u^2}-\frac{u_t}{u}\right)-\frac{n}{2}=-t\,\Delta\,\log u -\frac{n}{2}\leq 0\, .
\end{align}
Note, in particular, that the Li-Yau gradient estimate also implies that $F\leq 0$ when $\Ric\geq 0$ and hence $S\downarrow$.  
Li-Yau originally proved their gradient estimate using the maximum principle, but by now it is known how to deduce parabolic gradient estimates from monotonicity using the Shannon entropy and Fisher information and the monotonicity of $W$; see, for instance, \cite{BL}, \cite{N}, \cite{P2}.  As an almost immediate consequence of their gradient estimate, Li-Yau
 got a sharp Harnack inequality on manifolds with nonnegative Ricci curvature.  A key point in both the Li-Yau gradient estimate and the Harnack inequality is that they are sharp on Euclidean space for the heat kernel.  In fact, it is sharp in a very strong sense, namely if equality holds at one point on the manifold, then it is flat Euclidean space and $u$ is the heat kernel.

\subsection{Elliptic operators}
We will next turn to some elliptic analogs of the functionals and formulas defined in the previous subsection.  In the next section, we will give some applications of these.  This material is from \cite{C2}.   For simplicity, 
we will throughout the rest of the paper assume that $n\geq 3$.

Suppose that $u>0$ is harmonic in a pointed neighborhood of $p\in M$ and define $b$ by 
\begin{align}
	b^{2-n} =u \, . \notag 
\end{align}
  On Euclidean space,
 where $u=|x|^{2-n}$ is harmonic on $\RR^n\setminus \{0\}$,
  we get that $b=|x|$.   
The function $b$ satisfies the following crucial formula %that is sharp on the Green's function on Euclidean space
\begin{align} 
       \Delta \left(|\nabla b|^2\, u\right)=\frac{1}{2}\left(\left|\Hess_{b^2}-\frac{\Delta b^2}{n}g\right|^2+\Ric (\nabla b^2,\nabla b^2)\right)\,b^{-n}\, .  \label{e:bsub}
\end{align}
The right-hand side is always nonnegative on a smooth manifold $M$ with $\Ric \geq 0$ and it vanishes if and only if
$M$ is Euclidean space and $u$ is a multiple of  the Green's function.
We define a functional $A$ by
\begin{align}
A(r)&=r^{1-n}\int_{b=r}|\nabla b|^3\, .
\end{align}
A computation gives that the derivative  $A'(r)$ is equal to
\begin{align}
      - \frac{r^{n-3}}{2}\int_{b\geq r}\left(\left|\Hess_{b^2}-\frac{\Delta b^2}{n}g\right|^2+\Ric (\nabla b^2,\nabla b^2)\right)\,b^{2-2n}\, .		\notag
     \end{align}
     In particular, $A\downarrow$ is monotone nonincreasing on a manifold with $\Ric\geq 0$.  In the definition of $A$, we have implicitly assumed that $b$ is proper so that the integration is over a compact set.  This is automatically the case when $M$ is nonparabolic which is roughly equivalent to that the volume growth is faster than quadratic.  In the main application later on, $M$ will be assumed to have Euclidean volume growth and hence will be nonparabolic since $n\geq 3$. 

Often we will assume that $u$ is normalized so that it has the same asymptotics near the pole $p$ as the Green's function $|x|^{2-n}$ on Euclidean space near the origin.   If it is normalized in this way, then using \eqr{e:bsub} and the maximum principle,
 one can prove the sharp gradient estimate 
\begin{align}
|\nabla b|\leq 1 
\end{align}
on a manifold with $\Ric\geq 0$.  Similarly to the Li-Yau gradient estimate, if equality holds at one point on the manifold, then it is flat Euclidean space and $u$ is the Green's function.

 From the formula for $A'$ it follows, in particular, that if $M^n$ has $\Ric\geq 0$, then $A\downarrow$.   As $r$ tends to $0$, this quantity on a smooth manifold converges to the volume of the unit sphere in $\RR^n$ and, as $r$ tends to infinity, it converges to 
       \begin{align}
       \Vol (\partial B_1(0))\, \left(\frac{\V_M}{\Vol (B_1(0))}\right)^{\frac{2}{n-2}}\, .\label{e:asymptotic2}
       \end{align}
       Here, $\V_M$ is a geometric quantity that measures the Euclidean volume growth at infinity and will be defined in the next section and $B_1(0)\subset \RR^n$ is the unit ball in Euclidean space.
       In fact, one even has that
       \begin{align}
        \lim_{r\to \infty} \sup_{M\setminus B_r(x)}|\nabla b|=\left(\frac{\V_M}{\Vol (B_1(0))}\right)^{\frac{1}{n-2}}\, .
        \end{align}
        The function $b$ was already considered in \cite{CM1} where it was shown that on a manifold with $\Ric\geq 0$
       \begin{align}
        \lim_{r\to \infty} \sup_{\partial B_r(x)}\frac{b}{r}=\left(\frac{\V_M}{\Vol (B_1(0))}\right)^{\frac{1}{n-2}}\, .
        \end{align}

\section{Geometric inequalities and monotonicity}

We turn next to geometric applications of monotonicity and, in particular, how one can use the functionals and formulas of the previous section to prove uniqueness of blow-ups and blow-downs of Einstein manifolds.   Uniqueness is a key question for the regularity of Geometric PDE's.

\subsection{Scale invariant volume monotonicity and consequences}  The most basic monotonicity for Ricci curvature is that of the scale invariant volume.  This is usually called the Bishop-Gromov volume comparison theorem.  On its own it implies volume doubling and hence metric doubling and this gives directly compactness for the space of manifolds of a given dimension and a given lower bound of Ricci curvature.  To explain this,
 we need to recall a natural metric on the space of metric spaces.  This is the Gromov-Hausdorff distance that is a generalization of the classical Hausdorff distance between two subsets of the same Euclidean space.   Suppose that $(X,d_X)$ and $(Y,d_Y)$  are two compact metric spaces.  The Gromov-Hausdorff distance between them, denoted by $d_{GH} (X,Y)$, is by definition the infimum over all $\epsilon>0$ such that $X$ and $Y$ isometrically embed in a larger metric space $(Z,d_Z)$ and $X$ lies within an $\epsilon$-tubular neighborhood of $Y$ and vice versa.  A sequence of compact metric spaces $(X_i,d_{X_i})$ is said to converge in the Gromov-Hausdorff topology to a compact metric space $(Y,d_Y)$ if $d_{GH}(X_i,Y)\to 0$.  If the limit is noncompact, then the convergence is on compact subsets.

Gromov's compactness theorem is the result that any sequence of manifolds of a given dimension and a given lower Ricci curvature bound has a subsequence that converges in the Gromov-Hausdorff topology to a length space.  In particular, if $M$ has $\Ric\geq 0$, then any sequence of rescalings $(M,r_i^{-2}g)$, where $r_i\to \infty$, has a subsequence that converges in the Gromov-Hausdorff topology to a length space.   Any such limit is said to be a tangent cone at infinity of $M$.  As mentioned compactness follows from doubling that is implied by monotonicity of the scale invariant volume.  For a manifold with $\Ric\geq 0$, this is that
\begin{equation}
r^{-n}\,\Vol (B_r(x))
\end{equation}
is monotone nonincreasing in the radius $r$ of the ball $B_r(x)$ for any fixed $x\in M$.   As $r$ tends to $0$, this quantity on a smooth manifold converges to the volume of the unit ball in $\RR^n$ and, as $r$ tends to infinity, it converges to a nonnegative number $\V_M$.  If $\V_M>0$, then we say that $M$ has Euclidean volume growth and, by \cite{ChC1}, any tangent cone at infinity  is a metric cone.\footnote{A metric cone $C(X)$ with cross-section $X$ is a warped product metric $dr^2+r^2\,d^2_X$ on the space $(0,\infty)\times X$.  For tangent cones at infinity of manifolds with $\Ric\geq 0$ and $\V_M>0$, by \cite{ChC1}, any cross-section is a length space with diameter $\leq \pi$; see the next subsection for the precise definition of a metric cone over a general metric space.}

\subsection{$\Theta_r$ and geometric meaning of $A'$ and $W'$}
Our next goal is to explain the geometric meaning of $A'$ and similarly of $W'$ and for that we will need to recall what a metric cone is.   A metric cone $C(Y)$ over a metric space $(Y,d_Y)$ is the metric completion of the set $(0,\infty)\times Y$ with the metric
\begin{align}
d^2_{C(Y)}((r_1,y_1),(r_2,y_2))=r_1^2+r_2^2-2\,r_1\,r_2\,\cos d_Y(y_1,y_2)\, ;		\notag
\end{align}
see also section 1 of \cite{ChC1}.  When $Y$ itself is a complete metric space, taking the completion of $(0,\infty)\times Y$  adds only one point to the space.  This one point is usually referred to as the vertex of the cone.  We will also sometimes write $(0,\infty)\times_r Y$ for the metric cone.

We will next define the scale invariant distance between an annulus and the annulus in a cone centered at the vertex that approximates the annulus best.   (By scale invariant distance, we mean the distance 
   between the annuli after the metrics are rescaled so that the annuli have unit size.)  
To make this precise, suppose that $(X,d_X)$ is a metric space and $B_r(x)$ is a ball in $X$.  Let $\Theta_r(x)>0$ be the infimum of all $\Theta>0$ such that 
\begin{equation}
d_{GH}(B_{4r}(x)\setminus B_r(x),B_{4r}(v)\setminus B_r(v))<\Theta\,r\, ,	\notag
\end{equation}
where $B_r(v)\subset C(Y)$ and $v$ is the vertex of the cone.

Suppose now again that $u>0$ is harmonic in a pointed neighborhood of $p\in M$ and define $b$ by $b^{2-n} =u$ and $A$ by $A(r)=r^{1-n}\int_{b=r}|\nabla b|^3$.  
In the next subsection, when we discuss applications to uniqueness of blow-downs or blow-ups, it will be key that $A=A(r)$ is monotone $A\downarrow$ on a manifold with $\Ric\geq 0$ and for some positive constant $C$ 
 \begin{align}  \label{e:crucial}
 -A' (r)\geq C\,\frac{\Theta_r^2}{r} \, .
\end{align}
The constant $C$ depends only on a lower bound for the dimension and the scale invariant volume \footnote{The actual inequality is slightly more complicated as in reality the right hand side of this inequality is not to the power $2$ rather to the slightly worse power $2+2\epsilon$ for any $\epsilon> 0$, and the constant $C$ also depends on $\epsilon$; see \cite{C2}.}. This is a prime example of how the monotone quantities given purely analytically give information about the geometry of the spaces.  There is a similar inequality involving $W'$ and $\Theta_r$; see \cite{C2}.

\subsection{Uniqueness of tangent cones}
The quantity $A$ described above was used in \cite{CM2} to show that, for any Ricci-flat manifold with Euclidean volume growth, tangent cones at infinity are unique as long as one tangent cone has a smooth cross-section.  A similar result holds for local tangent cones of noncollapsed limits of Einstein manifolds.  Einstein manifolds and Ricci-flat manifolds, in particular, arise in a number of different fields, including string theory, general relativity, and complex and algebraic geometry, amongst others, and there is an extensive literature of examples.  Uniqueness of tangent cones is a key question for the regularity of Geometric PDE's.

There is a rich history of uniqueness results for geometric problems and equations. Ê In perhaps its simplest form, the issue of uniqueness   comes up already in a 1904 paper entitled ``On a continuous curve without tangents constructible from elementary geometry" by the Swedish mathematician Helge von Koch.  In that paper, Koch described what is now known as the Koch curve or Koch snowflake.  It is one of the earliest fractal curves to be described and, as suggested by the title, shows that there are continuous curves that do not have a tangent in any point.  On the other hand, when a set or a curve has a well-defined tangent or well-defined blow-up at every point, then much regularity is known to follow.  Tangents at every point, or uniqueness of blow-ups, is a `hard' analytical fact that most often is connected with a PDE,  as opposed to say Rademacher's theorem, where tangents are shown to exist almost everywhere for any Lipschitz function. 

In many geometric problems, existence of tangent cones comes from  monotonicity, while the approaches to uniqueness rely on showing that the monotone quantity approaches its limit at a definite rate.  However, estimating the rate of convergence seems to require either integrability and/or a great deal of regularity (such as analyticity).   For instance, for minimal surfaces or harmonic maps, the classical monotone quantities are highly regular and are well-suited to this type of argument.   However, this is not at all the case in the current setting where the Bishop-Gromov is of very low regularity and ill suited: the distance function is Lipschitz,  but is not even $C^1$, let alone analytic.  This is a major point.  In contrast, the functional $A$ is defined on the level sets of an analytic function (the Green's function) and does depend analytically and, furthermore, its derivative has the right properties.  In a sense, the scale invariant volume is already a regularization of $\Theta_r$ that, if one could, one would most of all like to work directly with and show some kind of decay for (in the scale).  However, not only is it not clear that $\Theta_r$ is monotone, but as a purely metric quantity it is even less regular than the scale invariant volume.
  
\subsection{Proving uniqueness}  To explain some of the key points of how one shows uniqueness of tangent cones, we let $p\in M$ be a fixed point in a Ricci-flat manifold with Euclidean volume growth.  We would like to show that the tangent cone at infinity is unique; that is, does not depend on the sequence of blow-downs.  To show this,
   let again $\Theta_r$ be the scale invariant Gromov-Hausdorff distance between the annulus $B_{4r}(p)\setminus B_{r}(p)$ and the corresponding annulus centered at the vertex of the cone that best approximates the annulus.   The first key point is that if $A=A(r)$ is defined as above, then $A$ is monotone, $A\downarrow$, and we have \eqr{e:crucial}.
(Perelman's monotone $W$ functional is also potentially a candidate, but it comes from integrating over the entire space which introduces so many other serious difficulties that it cannot be used.)    In fact, we shall use that, for $Q$ roughly equal to $-r\,A'(r)$, $Q$ is monotone nonincreasing and 
\begin{align}  \label{e:crucial2}
 [Q(r/2)-Q(8r)] \geq C\,\Theta_r^2 \, .
\end{align}
We claim that uniqueness of tangent cones is implied by showing that $A$ converges to its limit at infinity at a sufficiently fast rate or, equivalently, that $Q$ decays sufficiently fast to zero.  Namely, by the triangle inequality, uniqueness is implied by proving that
 \begin{align}
 \sum_k \Theta_{2^k}<\infty \, .
\end{align}
This, in turn, is implied by the Cauchy-Schwarz inequality by showing that for some $\epsilon>0$
\begin{align} \label{e:summa}
 \sum_k\Theta_{2^k}^2\,k^{1+\epsilon}<\infty \, ,
\end{align}
as
\begin{align} 
 \sum_k k^{-1-\epsilon}<\infty \, .
\end{align}
Equation \eqr{e:summa} follows, by \eqr{e:crucial2}, from showing that
\begin{align}  
 \sum [Q(2^{k-1})-Q(2^{k+3})]\,k^{1+\epsilon}<\infty \, .
\end{align}
This is implied by proving that, for a slightly larger $\epsilon$,
\begin{align}  
 Q(r)\leq  \frac{C}{(\log r)^{1+\epsilon}} \, .
\end{align}
All the work is then to establish  this crucial decay for $Q$.  This decay follows with rather elementary arguments
from showing that, for some $\alpha < 1 $,
 \begin{align}			\label{e:o3}
	    Q(2\,r)^{2-\alpha} &  \leq   C \, \left( Q(r/2) - Q(2\,r) \right) 
 \, .
\end{align}
The proof of this comes from an infinite dimensional Lojasiewicz-Simon inequality; see \cite{CM2} for details.

Finally, note that it is well-known that uniqueness may fail without the two-sided bound on the Ricci curvature.  Namely, there exist a large number of examples of manifolds with nonnegative Ricci curvature and Euclidean volume growth and nonunique tangent cones at infinity; see  \cite{ChC2}, \cite{CN}, \cite{P2}.

\end{document}